\documentclass[12pt]{amsart}
\begin{document}

\newcommand{\F}{\mathbb F}
\newcommand{\R}{\mathbb R}

\title[Logarithm of the Frobenius morphism]
{Heuristic formula for logarithm of the Frobenius morphism}
\author{A. Stoyanovsky}
\email{alexander.stoyanovsky@gmail.com}
\address{Russian State University of Humanities}

\begin{abstract}
We show that the logarithm $\log_q$ of the Frobenius morphism $x\to x^q$ is given by the formula
$x\to x\log x$ (the natural logarithm). In particular, it does not depend on $q$.
This is the explicit (although heuristical) formula for the
operator conjectured by Hilbert whose eigenvalues coincide with the zeroes of the zeta function.
\end{abstract}

\maketitle

The aim of this note is to make a remark clarifying the relation between the Riemann hypothesis on the
zeroes of the zeta function [1] and the Weil conjectures proved by P.~Deligne (see, for example, [2] for a brief
exposition). When trying to prove the Riemann hypothesis, it is natural to take into account its generalizations
to algebraic varieties, in particular to varieties over finite field $\F_q$. These Weil conjectures state that
the analog of the zeta function of a variety $X$ over $\F_q$, denote it by $\zeta_X(u)$, is a rational function of $u$
whose zeroes and poles lie on the circles of radius $|u|=q^{j/2}$, $0\le j\le\dim X$, and coincide with the eigenvalues
of the Frobenius morphism $F:x\mapsto x^q$ on the $j$-th $l$-adic cohomology of $X$. The analogy with the
Riemann hypothesis on the Riemann zeta function $\zeta(s)$ is achieved by putting $u=q^s$. Further, Hilbert
conjectured that there exists a Hilbert space $H$ with a self-adjoint operator $A$ with discrete spectrum in $H$
such that the eigenvalues of $A$ coincide, up to the transformation $s\to1/2+is$, with the critical zeroes of $\zeta(s)$.

Comparing these conjectures, one finds it natural to find the operator $B$ on cohomology of $X$ such that
$q^B=F$. The aim of this note is to write down a heuristic formula for the operator $B$. This can be fruitful for the
proof of the Riemann hypothesis.

Consider a 1-parametric discrete semigroup generated by $F$, $F^n:x\mapsto x^{q^n}$. Assume, heuristically,
that this 1-parametric family is extended to a continuous semigroup, $F^t:x\mapsto x^{q^t}$, $t\in\R$, $t\ge 0$. Let us
find the infinitesimal generator of this semigroup $B$ such that $F^t=\exp(Bt\log q)$. To this end, it suffices to
differentiate $F^t$ with respect to $t$:
\begin{equation}
\frac1{\log q}\frac d{dt}x^{q^t}=\frac1{\log q}x^{q^t}\log x\cdot q^t\log q=x^{q^t}\log x^{q^t}.
\end{equation}
In particular, for $t=0$ we obtain
\begin{equation}
B:x\mapsto x\log x.
\end{equation}
As expected, this formula does not depend on $q$.

The aim of further work would be to give a (homological) sense to formula (2). The operator $B$ should coincide
with the above mentioned operator $A$ (up to the transform $s\mapsto1/2+is$).

Note that this operator $B$ is almost a derivation, i.~e., it satisfies the property
\begin{equation}
B(xy)=Bx\cdot y+x\cdot By,
\end{equation}
but it is not linear: $B(x+y)\ne Bx+By$. This could be related to the fact that the difference $B(x+y)-Bx-By$ is
cohomologous to zero for any $x,y$.


\begin{thebibliography}{9}
\bibitem{1} E. C. Titchmarsh, The theory of the Riemann zeta function, Oxford, 1951.
\bibitem{2} R. Hartshorne, Algebraic geometry, Springer, 1977.
\end{thebibliography}
\end{document}